\documentclass{amsart}
\usepackage{amsmath,amssymb,amscd,amsthm}
\usepackage[matrix,arrow]{xy}

\newcommand{\g}{\mathfrak{g}}

\newcommand{\f}{\mathfrak{f}}
\newcommand{\fr}{\mathfrak{r}}
\newcommand{\nn}{\mathbb{N}}
\newcommand{\zz}{\mathbb{Z}}
\newcommand{\qq}{\mathbb{Q}}

\newcommand{\A}{\mathcal{A}}
\newcommand{\B}{\mathcal{B}}
\newcommand{\F}{\mathcal{F}}

\newtheorem{thm}{Theorem}
\newtheorem{pr}[thm]{Proposition}

\newtheorem{lm}{Lemma}[thm]

\theoremstyle{definition}
\newtheorem{df}{Definition}

\theoremstyle{remark}
\newtheorem*{rem}{Remark}

\newcommand{\la}{\lambda}
\newcommand{\ka}{\varkappa}
\newcommand{\ep}{\varepsilon}
\newcommand{\vphi}{\varphi}

\newcommand{\uu}{\mathbin{\smallsmile}}
\newcommand{\dd}{\partial}
\newcommand{\0}[1]{\overline{#1}}

\newcommand{\Hom}{\operatorname{Hom}}
\newcommand{\Cone}{\operatorname{Cone}}

\newcommand{\Ker}{\operatorname{Ker}}

\newcommand{\id}{\operatorname{id}}

\newcommand{\is}{\mathrel{\rlap{\raisebox{0.4ex}{\hspace{.4em}$\sim$}}%
\raisebox{-0.2ex}{$\longrightarrow$}}}

\begin{document}

\title{On the fundamental group and triple Massey's product}

\author{Grigori Rybnikov}
\thanks{Research
was supported in part by Grant M8H000/M8H300
from the International Science Foundation and Russian Government
and by INTAS Grant 94-4720.}
\address{Independent University of Moscow,
Department of Mathematics,
11 B. Vlas'evskij, Moscow 121002, Russia}

\email{gr@ium.ips.ras.ru}
\date{}

\maketitle

\section*{Introduction}

We study the relations between the fundamental group and the homological
operations on integer homology. For the ``rational'' fundamental group
(Malcev completion) see \cite{C1,C2,GM,Q,Su1,Su2}.

This work is an attempt to understand the invariant of the fundamental
group of the complement of a complex hyperplane arrangement that was
used in \cite{Ry}. Note that this invariant necessarily vanishes over
$\qq$ (see \cite{Ko}).

All homology and cohomology groups are with integer coefficients.
By $|X|$ we denote the geometric realization of a simplicial set $X$.

\section{Pseudo-isomorphisms and pseudo-ho\-me\-o\-mor\-phisms}

\begin{thm}\label{th5}
Suppose an arcwise connected topological space $U$
has the homotopy type of CW complex. Let $G=\pi_1(U,u)$, and let
$Y=B_.G$ be the nerve of $G$. Then there is a continuous map $U\to
BG$ inducing the natural isomorphism $\pi_1(U,u)=G\is\pi_1(BG,())$
{\rm(}hence, also an isomorphism $H_1(U)\to H_1(BG)${\rm)} and an
epimorphism $H_2(U)\to H_2(BG)$.
\end{thm}

\begin{proof}
Since $U$ has the homotopy type of a
CW complex, $U$ is homotopy equivalent to $|S(U)|$, where
$S(U)$ denotes the simplicial set of all singular simplices in $U$
\cite[Ch.\ 4]{LW}. Since $U$ is arcwise connected, $|S(U)|$ is homotopy
equivalent to $|S_u(U)|$, where $u\in U$ is an arbitrary point and
$S_u(U)$ denotes the simplicial set of all singular simplices in $U$
with all vertices equal to $u$. Therefore, it suffices to prove the
theorem for $U=|X|$, where $X$ is a simplicial set with single vertex
$x$, $u=x$. 

In fact, we will construct a simplicial map $X\to
B_.G$ inducing the natural isomorphism $\pi_1(|X|,x)=G\is\pi_1(BG,())$
and an epimorphism $H_2(|X|)\to H_2(BG)$.

Let $\tilde X_k$ denote the set of non-degenerate $k$-simplices in $X$. 
Consider the free group $F=F(\tilde X_1)$ and its nerve $B_.F$. It is
well known that $H_0(BF)=\zz$, $H_1(BF)=\zz^{\tilde X_1}$, and
$H_i(BF)=0$ for $i>1$. Denote by $X^1$ the minimal simplicial subset of
$X$ containing $\tilde X_1$. We identify $X^1$
with the corresponding simplicial subset of $B_.F$ and glue $|X|$ and
$BF$ via this identification. Denote the resulting simplicial set by
$Y$.

It is clear that $|X^1|$ is a union of
1-dimensional spheres, hence the inclusion $|X^1|\hookrightarrow BF$
induces isomorphism $H_i(|X_1|)\is H_i(BF)$ for any $i$. Hence, from
the long
homological exact
sequence of the pair $(BF,|X^1|)$ we
see that $H_i(BF,|X^1|)=0$ for any $i$. But $H_i(|Y|,|X|)=H_i(BF,|X^1|)$.
From the long
homological exact
sequence of the pair $(|Y|,|X|)$ we
see that the natural map $H_i(|X|)\to H_i(|Y|)$ is an isomorphism.
It is also clear that the inclusion of $X$ into $Y$ gives an isomorphism of
fundamental groups. Thus we can use $Y$ instead of $X$.

Let us define a simplicial map $\varphi Y\to B_.G$. We have a natural
homomorphism $F\to G$ that maps each element $a\in\tilde X_1$ to the
corresponding element $g_a$ of the fundamental group $\pi_1(X,x)=G$. We
put $\varphi((a_1,\dots,a_n))=(g_{a_1},\dots,g_{a_n})$ for each
$(a_1,\dots,a_n)\in B_nF$. Suppose $\sigma$ is a $k$-simplex in $X$. Let
$a_k$ be $1$-face of $\sigma$ corresponding to the inclusion
$f_k:[0,1]\to{}[0,m]$ given by $f_k(j)=j+k$. We put
$\varphi(\sigma)=(g_{a_1},\dots,g_{a_k})$. It is obvious that $\varphi$
is a simplicial map. Besides, it is clear that the map
$\varphi_*:\pi_1(|Y|,x)\to\pi_1(BG,())$ is an isomorphism. 

Note that the map $\varphi_*:C(Y)\to C(B_.G)$ is surjective.
Denote by $K$ the kernel of this map. We have $K_0=0$ and
$K_1=\langle (gn)-(g)\rangle_{n\in N}$, where $N$ is the kernel of the
natural homomorphism $F\to G$.
From the long exact sequence
$$
\dots\to H_2(Y)\to H_2(BG)\to H_1(K)\to\dots
$$
we see that it suffices to show that
$H_1(K)=0$.

It is clear that for any $g_1,g_2,h\in F$ such that $(g_1)-(g_2)\in\dd
K_2$ we have $(hg_1)-(hg_2)\in\dd K_2$ and $(g_1h)-(g_2h)\in\dd K_2)$.
It follows that for any $g,h\in F$ such that $(g)-(1)\in\dd K_2$ we have
$(h^{-1}gh)-(1)\in\dd K_2$ and for any $g,h\in F$ such that
$(g)-(1)\in\dd K_2$ and $(h)-(1)\in\dd K_2$ we have
$(gh)-(1)=(gh)-(h)+(h)-(1)\in\dd K_2$.

Let us consider arbitrary element $\sigma\in X_2$ and
let $a=d_2\sigma$, $b=d_0\sigma$, and $c=d_1\sigma$. Denote the
corresponding elements of $F$ by $g_a$, $g_b$, and $g_c$. We have
$(g_ag_b)-(g_c)=\dd((g_a,g_b)-\sigma)\in \dd K_2$. Hence
$(g_ag_bg_c^{-1})-(1)\in\dd K_2$. The elements of the form
$g_ag_bg_c^{-1}$ and their conjugates generate the subgroup $N$.
Therefore for any $n\in N$ we have $(n)-(1)\in\dd K_2$. It follows that 
$(gn)-(g)\in\dd K_2$ for any $n\in N$, $g\in F$, thus $K_1=\dd K_2$.
Hence $H_1(K)=0$ and the map $\varphi_*:H_2(|Y|)\to H_2(BG)$ is an
epimorphism.
\end{proof}

Let us consider only arcwise connected topological spaces that
have the homotopy type of CW complex. Any continuous map
inducing isomorphism of $H_1$ and epimorphism of $H_2$ will
be called {\em a pseudo-homeomorphism\/}. We see that any invariant of
topological space w.r.t.\ pseudo-homeomorphisms is an invariant of its
fundamental group.

Now we want to know what information about the fundamental group
can be contained in such invariants.

\begin{df}
 Let $G$ be a group. Denote by $I$ the augmentation ideal of the group
algebra $\zz G$ (that is, $I$ is generated by all elements of the form
$g-e$, where $g\in G$ and $e$ is the identity element of $G$). We put
$D^{(k)}(G)=\zz G/I^k$.
Denote by $D(G)$ the projective system of $\zz$-algebras
$$\to D^{({k+1})}(G)\to D^{(k)}(G)\to D^{({k-1})}(G)\to\dots\to
D^{(1)}(G)=\zz.$$

Suppose $\varphi:G_1\to G_2$ is a group homomorphism. We say that it is
{\em a pseudo-isomorphism\/} if it gives rise to an isomorphism
$D(G_1)\to D(G_2)$.
\end{df}

Let $X=(X_n)_{n\in\zz_+}$ be a simplicial set with single vertex (that is,
$X_0=\{x\}$). The
fundamental group of $|X|$ can be described as the group with generators
$g_a$ ($a\in X_1$) and relations $g_{d_2\sigma}g_{d_0\sigma}=
g_{d_1\sigma}$ for any $\sigma\in X_2$, where $d_i\sigma$ means $i$-th
face of $\sigma$ (with $i$-th vertex missing).

Let $C=C(X)$ be the chain complex of $X$ over $\zz$. It has the standard
structure of coalgebra: for any $n$-simplex $\sigma$ we have
$$
\Delta\sigma=\sum_{i+j=n}{}_{(i)}\sigma\otimes\sigma_{(j)},
$$
where ${}_{(i)}\sigma$ is the front $i$-dimensional face and
$\sigma_{(j)}$ is the back $j$-dimensional face of $\sigma$.
Let $\0C=C/C_0$; since $C_0=\zz x$ is subcoalgebra of $C$, we obtain
comultiplication $\0\Delta:\0C\to\0C\otimes\0C$. Denote by
$\F(C)$ the tensor algebra $T(s^{-1}\0C)$ (cobar construction
\cite{A}). Note that $\F(C)_0$ is a free associative algebra generated by
$X_1$. We write $[c_1|c_2|\dots|c_k]$ instead of
$s^{-1}c_1\otimes s^{-1}c_2\otimes\dots\otimes s^{-1}c_k$.
The differential in $\F(C)$ is a derivation of the tensor algebra
defined on generators as
$$
\dd [c]=-[\delta c]+\sum(-1)^{\deg a_i}[a_i|b_i],
$$
where $\0\Delta c=\sum a_i\otimes b_i$ and $\delta$ is the differential
in $\0C$.

Let $\F^{(k)}(C)=\F(C)/(T^k(s^{-1}\0C))$. It is clear that the ideal
$(T^k(s^{-1}\0C))=\bigoplus_{r\ge k}T^r(s^{-1}\0C)$ is a
subcomplex of $\F(C)$; hence $\F^{(k)}(C)$ is a complex. Denote
$H_0(\F^{(k)}(C))$ by $A^{(k)}=A^{(k)}(X)$.

\begin{pr}\label{digr}
Let $G=\pi_1(|X|,v)$. Then there is an isomorphism $D^{(k)}(G)\to
A^{(k)}(X)$ such that for each $a\in X_1$ the image of $g_a$ in
$D^{(k)}(G)$ corresponds to the image of $1+[a]$ in $A^{(k)}$.
\end{pr}

\begin{proof}
Denote $1+[a]\in \F(C)_0$ by $\tilde g_a$. Let $\sigma$ be a $2$-simplex and
let $a=d_2\sigma$, $b=d_0\sigma$, and $c=d_1\sigma$. Then $\tilde
g_a\tilde g_b-\tilde g_c=[a]+[b]+[a|b]-[c]=-\dd[\sigma]$, thus the
corresponding element in $A^{(k)}$ is zero. Hence there is a
homomorphism $\zz G\to A^{(k)}$ sending each generator $g_a$ to the
image of $\tilde g_a$ in $A^{(k)}$. From the definitions it is clear
that its kernel is $I^k$.
\end{proof}

\begin{thm}\label{th4}
Let $U$ and $V$ be arcwise connected topological spaces
having the homotopy types of CW complexes.
Suppose $f:U\to V$ is a pseudo-ho\-me\-o\-mor\-phism. Then
$$f_*:\pi_1(U,u)\to\pi_1(V,f(u))$$
is a pseudo-isomorphism.
\end{thm}

\begin{proof}
As in the proof of Theorem \ref{th5}, we assume that $U=|X|$ and
$V=|Y|$, where $X$ and $Y$ are simplicial spaces with single vertices
$x$ and $y$ respectively. Besides, we assume that $f=|F|$, where
$F:X\to Y$ is a simplicial map, $u=x$ and, thus, $f(u)=y$.

By Proposition \ref{digr} it suffices to show that the natural map
$A^{(k)}(X)\to A^{(k)}(Y)$ is an isomorphism for any $k\in\nn$.
Consider the filtration of the complex $\F^{(k)}(C)$ for both $C=C(X)$
and $C=C(Y)$
$$
\F^{(k)}(C)=\F_0^{(k)}(C)\supset \F_{-1}^{(k)}(C)\supset\dots
\supset \F_{-k+1}^{(k)}(C)\supset \F_{-k}^{(k)}(C)=\{0\}
$$
where $\F_p^{(k)}(C)=(T^{-p}(s^{-1}\0C)/(T^k(s^{-1}\0C)$. We have a
natural map of the corresponding spectral sequences $E_{p,q}(X)\to
E_{p,q}(Y)$.

Note that $E^1_{-r,q}=H_{q-r}((s^{-1}\0C)^{{\otimes}r})$.

\begin{lm}
Let $K$, $L$, and $M$ be chain complexes of free $\zz$-modules
and let $f:K\to L$ be a map of
complexes such that $f_*:H_i(K)\to H_i(L)$ is an isomorphism for
$i=1,\dots,k$ and an epimorphism for $i=k+1$. Then
$(f\otimes\id_M)_*:H_i(K\otimes M)\to H_i(L\otimes M)$ is an isomorphism for
$i=1,\dots,k$ and an epimorphism for $i=k+1$.
\end{lm}

\begin{proof}
The condition of the Lemma is equivalent to the fact that $\Cone(f)$ is
acyclic in dimensions $0,1,\dots,k$. But
$\Cone(f\otimes\id_M)\simeq\Cone(f)\otimes M$, therefore it is also
acyclic in dimensions $0,1,\dots,k$.
\end{proof}

Let us continue the proof of the Theorem.
Repeatedly applying the Lemma we see that the natural map
$E^1_{p,q}(X)\to E^1_{p,q}(Y)$ is an isomorphism for $p+q=0$ and an
epimorphism for $p+q=1$. It follows that it is also true for $E^r_{p,q}$
for any $r>0$, therefore it is true for $E^\infty_{p,q}$ and hence for
$H_i(\F^{(k)}(C))$.
\end{proof}

\section{Triple Massey's product}

From theorem \ref{th4} it follows that the invariants of a topological space
w.r.t.\ pseudo-homeomorphisms distinguish
fundamental groups (at least) up to a pseudo-iso\-mor\-phism.

Clearly, in our study of such invariants
it suffices to consider simplicial sets with single vertex.
A simplicial map $f:X\to Y$ of simplicial sets with single vertex will
be called a {\em simplicial pseudo-homeomorphism\/}, if its geometric
realization $|f|:|X|\to |Y|$ is a pseudo-homeomorphism.

\begin{pr}
Let $f:X\to Y$ be a simplicial pseudo-homeomorphism.
Then the map $f^*:H^i(Y)\to
H^i(X)$ is an isomorphism for $i=1$ and a monomorphism for $i=2$.
\end{pr}

\begin{proof}
Let us consider the map $\vphi:C_.(X)\to C_.(Y)$ of chain complexes.
We put $C_.=C_.(X)$, $K_.=C_.(Y)$. Thus
$C_0\simeq\zz$, $K_0\simeq\zz$, $\dd_{C_1}=0$,  $\dd_{K_1}=0$. 
We have
\begin{equation}\label{cond}
K_1=\vphi C_1+\dd K_2,\quad\vphi^{-1}(\dd K_2)=\dd
C_2,\quad\Ker\dd_{K_2}=\vphi\Ker\dd_{C_2}+\dd K_3.
\end{equation}

Let $C^{{\cdot}}=\Hom(C_.,\zz)$ and $K^{{\cdot}}=\Hom(K_.,\zz)$. From the
first equation in \eqref{cond} we get $\Ker\vphi^*\cap\Ker d_{K^1}=0$,
combining the first and the second equations in \eqref{cond} we get
$\Ker d_{C^1}=\vphi^*\Ker d_{K^1}$, and with the help of all
\eqref{cond} we get ${\vphi^*}^{-1}(dC^1)\cap\Ker d_{K^2}=dK^1$.
This is just what we need.
\end{proof}

Let $X$ be a simplicial set.
We write $H^i$ instead
of $H^i(X)$ and use standard notations $C^i$, $Z^i$, and $B^i$ for
$\zz$-modules of $i$-cochains, $i$-cocycles, and $i$-coboundaries
respectively. Let us denote by $\mu$ the map of cup-product
$C^i\otimes C^j\to C^{i+j}$ (the map $\mu:C^{\cdot}\otimes C^{\cdot}\to
C^{\cdot}$ is the adjoint operator to comultiplication $\Delta:C_.\to
C_.\otimes C_.$). The cup-product in cohomology $H^i\otimes H^j\to
H^{i+j}$ is denoted by $\bar\mu$.

Now we will construct an invariant of pseudo-homeomorphisms that is
closely related to the triple Massey product.
Let us recall its definition.
Suppose $\theta_1,\theta_2$, $\theta_3\in H^1$ satisfy the conditions $\theta_1\uu \theta_2=0$ and
$\theta_2\uu \theta_3=0$.
Choose $\omega_i\in Z^1$ such that
$\theta_i=[\omega_i]$ ($i=1,2,3$); then there exist $\omega_{12},\omega_{23}\in
C^1$ such that $\omega_1\uu\omega_2=d\omega_{12}$ and
$\omega_2\uu\omega_3=d\omega_{23}$. Consider the cochain
$\omega_{12}\uu\omega_3+\omega_1\uu\omega_{23}$; clearly, it is a cocycle; its
class in $H^2$ is denoted by $\langle \theta_1,\theta_2,\theta_3\rangle$ and is
called {\em triple Massey product\/} of $\theta_1,\theta_2,\theta_3$. This product is
not defined uniquely: only the set $\langle \theta_1,\theta_2,\theta_3\rangle+H^1\uu
\theta_3+\theta_1\uu H^1$ has an invariant sense.

Denote by $\eta_2$ the natural projection $H^1\otimes
H^1\to\Lambda^2H^1$; choose a homomorphism of $\zz$-modules
$\chi_2:\Lambda^2H^1\to H^1\otimes H^1$ such that
$\eta_2\circ\chi_2=\id$ (in other words, $\chi_2$ is right inverse to
$\eta_2$). Since $H^1$ is a free Abelian group,
we can fix a homomorphism $\ka:H^1\to Z^1$ right inverse to the
canonical projection.

\begin{pr}\label{zeta}
There is a natural (non-linear) map $\zeta:Z^1\to C^1$
such that $\omega\uu\omega=d\zeta(\omega)$ for any $\omega\in
Z^1$.
\end{pr}

\begin{proof}
The map $\zeta$ may be defined as follows:
$\zeta(\omega)(a)=(\omega(a)-\omega(a)^2)/2$ for any 1-simplex $a$. It
is clear that this map is natural; the equality
$\omega\uu\omega=d\zeta(\omega)$ is easy to check.
\end{proof}

From Proposition \ref{zeta} it follows that the map $\bar\mu:H^1\otimes
H^1\to H^2$ may be factored through $\Lambda^2H^1$; it means that there
is a map $\Bar{\Bar\mu}:\Lambda^2H^1\to H^2$ such that
$\Bar{\Bar\mu}\circ\eta_2=\bar\mu$. Let $R^2=\Ker{\Bar\mu}$,
$\bar R^2=\Ker\Bar{\Bar\mu}$. Clearly, $S^{(2)}H^1\subset R^2$, where
$S^{(2)}H^1\subset H^1\otimes H^1$ is the set of symmetric tensors.

Denote $R^2\otimes H^1\cap H^1\otimes R^2\subset H^1\otimes H^1\otimes
H^1$ by $Q^3$. Let us choose a homomorphism of Abelian groups
$\nu:R^2\to C^1$ such that $\mu\circ(\ka\otimes\ka)(r)=d\nu(r)$ for any
$r\in R^2$. Now define $\lambda:Q^3\to C^2$ as
$\lambda=\mu\circ(\nu\otimes\ka+\ka\otimes\nu)$.
We have
$$
d\circ\lambda=\mu\circ((d\circ\nu)\otimes\ka-\ka\otimes(d\circ\nu))
=\mu\circ(\mu\circ(\ka\otimes\ka)\otimes\ka
-\ka\otimes\mu\circ(\ka\otimes\ka))=0
$$
by associativity of the cup-product. For any $q\in Q^3$ we denote the
image of $\lambda(q)$ in $H^2$ by $\bar\lambda(q)$.

Let $\theta_1$, $\theta_2$, and $\theta_3$ be as in the definition of the triple Massey
product. Then $q=\theta_1\otimes \theta_2\otimes \theta_3\in Q^3$ and $\langle
\theta_1,\theta_2,\theta_3\rangle=\bar\lambda(q)$. Thus, the homomorphism
$\bar\lambda:Q^3\to H^2$ can be viewed as a form of the triple Massey
product.

Clearly, $\bar\lambda$ depends on the choices of $\ka$ and $\nu$.
If we change $\ka$ to $\ka'=\ka+d\circ\ep$, where $\ep:H^1\to C^0$ is an
arbitrary homomorphism, then
$\nu'=\nu+\mu\circ(\ka\otimes\ep)+\mu\circ(\ep\otimes\ka')$ satisfies
the condition $\mu\circ(\ka'\otimes\ka')=d\circ\nu'(r)$ and gives rise
to the same $\bar\lambda$. On the other hand, if we change $\nu$ to
$\nu'=\nu+\rho$, where $\rho:R^2\to Z^1$ is an arbitrary homomorphism,
then $\bar\lambda$ changes to
$\bar\lambda+\bar\mu\circ(\id\otimes\bar\rho+\bar\rho\otimes\id)$, where
$\bar\rho$ is the composition of $\rho$ with the canonical projection
$Z^1\to H^1$. Denote by $\delta$ the map from $\Hom(R^2,H^1)$ to
$\Hom(Q^3,H^2)$ sending $f\in\Hom(R^2,H^1)$ to
$\bar\mu\circ(\id\otimes f+f\otimes\id)$. We see that the class of
$\bar\lambda$ in $\Hom(Q^3,H^2)/\delta\Hom(R^2,H^1)$ is well-defined.
Clearly, it is an invariant of pseudo-homeomorphisms.

We can further reduce $\bar\lambda$ with the help of Proposition
\ref{zeta}.

Let us fix a basis $(\xi_1,\dots,\xi_n)$ of $H^1$. We choose $\nu$ so
that $\nu(\xi_i\otimes \xi_i)=\zeta(\ka(\xi_i))$ for $i=1,\dots,n$ and
$\nu(\xi_i\otimes \xi_j+\xi_j\otimes
\xi_i)=\zeta(\ka(\xi_i+\xi_j))-\zeta(\ka(\xi_i))-\zeta(\ka(\xi_j))$ for
$i\ne j$. Let $\bar\delta$ be the map from $\Hom(\bar R^2,H^1)$ to
$\Hom(Q^3,H^2)$ sending $f\in\Hom(\bar R^2,H^1)$ to
$\delta(f\circ\eta_2)$. Clearly, the class of $\bar\lambda$ in
$\Hom(Q^3,H^2)/\bar\delta\Hom(\bar R^2,H^1)$ is a well-defined invariant
of pseudo-homeomorphisms.

Now consider the map $l:H^1\otimes\bar R^2\to\Lambda^3H^1$ arising from
the wedge product in $\Lambda^{\bullet}H^1$.
We set $\bar Q^3=\Ker l$.

\begin{pr}
The image of the map $(\id_{H^1}\otimes\eta_2)\circ(\id-s_{(123)}):Q^3\to
H^1\otimes\bar R^2$ is $\bar Q^3$.
\end{pr}

\begin{proof}
Clearly, the image of this map belongs to $\bar Q^3$. Let us construct
a map $p:\bar Q^3\to Q^3$ right inverse to the map under consideration.

Let $t=\sum_i\xi_i\otimes r_i\in\bar Q^3$,
$r_i=\sum_{j<k}\alpha_{ijk}\xi_j\wedge \xi_k$. We put
\begin{multline*}
pt=\sum_{i<j<k}(\alpha_{jik}(\xi_i\otimes \xi_j\otimes \xi_k+\xi_j\otimes
\xi_i\otimes \xi_k)+\alpha_{kij}(\xi_i\otimes \xi_k\otimes \xi_j+\xi_k\otimes
\xi_i\otimes \xi_j))\\+\sum_{i<j}(\alpha_{iij}\xi_i\otimes \xi_i\otimes \xi_j
-\alpha_{jij}(\xi_j\otimes \xi_j\otimes \xi_i).
\end{multline*}
Since $\alpha_{ijk}-\alpha_{jik}+\alpha_{kij}=0$ for $i<j<k$, we have
$pt\in H^1\otimes R^2\cap S^{(2)}H_1\otimes H^1$ and  
$(\id_{H^1}\otimes\eta_2)\circ(\id-s_{(123)})pt=t$.
\end{proof}

Let us fix the map $p$ constructed above.

Let now $t\in H^1\otimes\bar R^2$, $t=\sum_{i}\xi_i\otimes r_i$,
$r_i=\sum_{j<k}\alpha_{ijk}\xi_j\wedge \xi_k$. We put
$$
qt=\sum_i\sum_{j<k}\alpha_{ijk}(\xi_i\otimes \xi_j\otimes \xi_k+\xi_j\otimes
\xi_i\otimes \xi_k+\xi_j\otimes \xi_k\otimes \xi_i).
$$
Clearly, $q$ is the map $H^1\otimes\bar R^2\to Q^3$ satisfying the
conditions $(\id_{H^1}\otimes\eta_2)\circ q=\id$ and
$(\id_{H^1}\otimes\eta_2)\circ s_{(123)}\circ q=\id$.

\begin{pr}
We have $Q^3=p\bar Q^3\oplus q(H^1\otimes\bar R^2)\oplus S^{(3)}H^1$.
\end{pr}

\begin{proof}
Clear.
\end{proof}

Now let $\bar\lambda:Q^3\to H^2$ be as above.

\begin{pr}
We have $\bar\lambda(t)=0$ for any $t\in S^{(3)}H^1$ and
$\bar\lambda(qt)=\sum_{i<j}(\alpha_{iij}+\alpha_{jij})\xi_i\uu \xi_j$ for
$t=\sum_{i}\xi_i\otimes\sum_{j<k}\alpha_{ijk}\xi_j\wedge \xi_k\in
H^1\otimes\bar R^2$.
\end{pr}

\begin{proof}
Let $h\in H^1$, $\tilde h=\ka h$. We define $f\in C^1$ by the formula
$f(a)=\tilde h(a)(\tilde h(a)-1)(\tilde h(a)-2)/6$ for all $a\in X_1$.
It is easy to check that $\lambda(h\otimes h\otimes h)=df$.
Therefore, $\bar\la$ vanishes on $S^{(3)}H^1$.

To prove the second assertion of the proposition we recall that there is
a natural map $\mu_1:C^{{\cdot}}\otimes C^{{\cdot}}\to C^{{\cdot}}$ of
degree $-1$ such that for any $f\in C^{{\cdot}}\otimes C^{{\cdot}}$ one
has $\mu(f-sf)=d\mu_1f+\mu_1df$, where $s(a\otimes b)=(-1)^{\deg a\deg
b}b\otimes a$ for homogeneous $a,b\in C^{{\cdot}}$. (This is a part of
the structure of $E_\infty$-algebra on $C^{{\cdot}}$, see \cite{S}.) The
map $\mu_1$ is adjoint to the map $\Delta_1:C_.\to C_.\otimes C_.$ of
degree $1$. For the standard simplices of dimensions $1$ and $2$ the map
$\Delta_1$ is given by the formulas $\Delta_1[01]=[01]\otimes[01]$ and
$\Delta_1[012]=[012]\otimes[02]+([01]+[12])\otimes[012]$.

For
$t=\sum_{i}\xi_i\otimes\sum_{j<k}\alpha_{ijk}\xi_j\wedge \xi_k\in
H^1\otimes\bar R^2$ we have
\begin{multline*}
\lambda(qt)=\mu\circ(\nu\otimes\ka-s\circ(\nu\otimes\ka)+(\ka\otimes\nu)
\circ(\id-s_{(123)}))(qt)\\
=d\circ\mu_1\circ(\nu\otimes\ka)(qt)+
\mu_1\circ(\mu\otimes\id)\circ(\ka\otimes\ka\otimes\ka)(qt)+(\ka\otimes\nu)
\circ(\id-s_{(123)})(qt).
\end{multline*}
Note that $(\id-s_{(123)})(qt)$ is symmetric w.r.t.\ the last two
indices. Now it is easy to check that
\begin{multline*}
\mu_1\circ(\mu\otimes\id)\circ(\ka\otimes\ka\otimes\ka)(qt)+(\ka\otimes\nu)
\circ(\id-s_{(123)})(qt)\\=-d\sum_{i}\sum_{j<k}\alpha_{ijk}\tilde
\xi_i\tilde \xi_j\tilde \xi_k+\sum_{i<j}(-\alpha_{iij}\tilde \xi_j\uu\tilde
\xi_i+\alpha_{jij}\tilde \xi_i\uu\tilde \xi_j),
\end{multline*}
where $\tilde h=\ka h$ and $\tilde
\xi_i\tilde \xi_j\tilde \xi_k$ is a cochain in $C^1$ given by $$\tilde
\xi_i\tilde \xi_j\tilde \xi_k(a)=\tilde
\xi_i(a)\tilde \xi_j(a)\tilde \xi_k(a)$$
for all $a\in X_1$.
\end{proof}

We see that $\bar\lambda$ is determined by $\Bar{\Bar\la}=\bar\la\circ
p:\bar Q^3\to H^2$. Denote by $\Bar{\Bar\delta}$ the map $\Hom(\bar
R^2,H^1)\to\Hom(\bar Q^3,H_2)$ sending $f\in\Hom(\bar R^2,H^1)$ to
$\bar\delta(f)\circ p$. Clearly, the class $[\la]\in\Hom(\bar
Q^3,H^2)/\Bar{\Bar\delta}\Hom(\bar R^2,H^1)$ of $\Bar{\Bar\la}$ is a
well-defined invariant of pseudo-homeomorphisms. More precisely, we have

\begin{thm}\label{th8}
Let $X$ and $Y$ be simplicial sets and let
$\vphi:X\to Y$ be simplicial pseudo-homeomorphism. Suppose that the map
$\chi_2^{(X)}:\Lambda^2H^1(X)\to H^1(X)\otimes H^1(X)$ and the basis
of $H^1(X)$ are chosen as above.
Let us transfer these
structures to $H^1(Y)$ with the help of isomorphism $\vphi^*:H^1(Y)\is
H^1(X)$ and construct $\bar R^2(X)$, $\bar R^2(Y)$,
$\bar Q^3(X)$, $\bar Q^3(Y)$, $[\la]^{(X)}$, $[\la]^{(Y)}$. Then
\begin{enumerate}
\item
the map $\Lambda^2\vphi^*:\Lambda^2H^1(Y)\to\Lambda^2H^1(X)$ restricted
to $\bar R^2(Y)$ is an isomorphism of $\bar R^2(Y)$ onto $\bar R^2(X)$;
\item
the map $\Lambda^2\vphi^*\otimes\vphi^*:\Lambda^2H^1(Y)\otimes
H^1(Y)\to\Lambda^2H^1(X)\otimes H^1(X)$ restricted to $\bar Q^3(Y)$
is an isomorphism of $\bar Q^3(Y)$ onto $\bar Q^3(X)$;
\item
$[\la]^{(X)}$ coincides with the image of $[\la]^{(Y)}$
under the map
$$
\Hom_\zz(\bar Q^3(Y),H^2(Y))\to\Hom_\zz(\bar Q^3(X),H^2(X))
$$
given by $\Bar{\Bar\lambda}\mapsto\vphi^*\circ\Bar{\Bar\lambda}\circ
(\Lambda^2\vphi^*\otimes\vphi^*)^{-1}$.
\end{enumerate}
\end{thm}

\begin{proof}
This is obvious since we have used only natural constructions to
define $\bar R^2$, $\bar Q^3$, and $[\lambda]$.
\end{proof}

\section{The invariant of fundamental group}

Our next goal is to interpret the invariant $[\lambda]$ in terms of
the fundamental group. It is possible due to theorem \ref{th5},
but we want to produce an explicit construction.

Let $G=\pi_1(U,u)$, where $U$ is an arcwise connected topological space
having the homotopy type of a CW-complex.
Suppose the following conditions hold:
\begin{enumerate}
\item
$G/G'$ is a free Abelian group of rank $n$ (that is,
$H_1(U)\simeq\zz^n$);
\item
$G$ is generated with $n$ generators, thus $G\simeq F/R$,
where $F=F(w_1,\dots,w_n)$ is a free group with generators $w_1,\dots,w_n$
and $R$ is a normal subgroup of $F$;
\item
the comultiplication in $C_.(U)$ gives rise to an injective homomorphism
$$H_2(U)\to H_1(U)\otimes H_1(U).$$
\end{enumerate}

\begin{rem}
These conditions are satisfied for the complement of a complex
hyperplane arrangement (see \cite{OT}).
\end{rem}

By Theorem \ref{th5} there is a pseudo-homeomorphism $U\to BG$. Hence, by
Theorem \ref{th8} the invariants $\bar R^2$, $\bar Q^3$, and $[\la]$ may be
computed for the simplicial set $X=B_.G$ with single vertex $x=()$.

Denote by $g_1,\dots,g_n$ the generators of $G$ (i.~e., the images of
$w_1,\dots,w_n$). Let $h_1\dots,h_n$ be the corresponding elements in
$H_1=H_1(X)\simeq G/G'\simeq F/F'$. Clearly, $(h_1,\dots,h_n)$ is a
basis of $H_1$. Since $H_2=H_2(X)$ is imbedded into $H_1\otimes H_1$, it
is also a free $\zz$-module. Hence $H_1\simeq\Hom_\zz(H^1,\zz)$ and
$H_2\simeq\Hom_\zz(H^2,\zz)$. Denote by $(\xi_1,\dots,\xi_n)$ the basis
of $H^1$ dual to $(h_1,\dots,h_n)$.

As usual, denote by $\gamma_kG$ the $k$-th term of the lower central
series of $G$ (that is, $\gamma_1G=G$ and
$\gamma_{k+1}G=(G,\gamma_kG)=\langle g^{-1}f^{-1}gf\mid g\in
G,f\in\gamma_kG\rangle$). It is well known that
$\g=\bigoplus_{k=1}^\infty\g_k$, where $\g_k=\gamma_kG/\gamma_{k+1}G$, has
the structure of the graded Lie algebra with the Lie commutator $[\ ,\
]:\g_k\times\g_m\to\g_{k+m}$ in $\g$ corresponding to the group
commutator $(\ ,\ ):\gamma_kG\times\gamma_mG\to\gamma_{k+m}G$.
Note also that for the free group $F$ the corresponding Lie algebra $\f$
is the free Lie algebra with generators $x_i=w_i\gamma_2F\in\f_1$ (see
\cite{MKS}). 

By the Magnus Theorem \cite{M,Sa} the subgroup $\gamma_kF$ is the set of
all $w\in F$ such that $w-1$ belongs to the $k$-th power of the
augmentation ideal $J\in\zz F$. This is not generally true for an
arbitrary group (see \cite{R}). But in our case we have the following 

\begin{pr}\label{gamma}
For $k=1,2,3,4$ the subgroup $\gamma_kG$ is the set of
all $g\in G$ such that $g-1$ belongs to the $k$-th power of the
augmentation ideal $I\in\zz G$.
\end{pr}

\begin{proof}
For $k=1$ there is nothing to prove. The case $k=2$ is not much harder.

Since $G=F/R$, we have $\gamma_kG=\gamma_kF/R\cap\gamma_kF$; thus,
$\g=\f/\fr$, where $\fr=\bigoplus\fr_k$,
$\fr_k=R\cap\gamma_kF/R\cap\gamma_{k+1}F$.

Consider the graded algebra $\A=\sum_{k=0}^\infty J^k/J^{k+1}$. Clearly,
it is a free associative algebra with the generators $x_1,\dots,x_n$,
where $x_i=(w_i-1)+J_2$; $\A$ is isomorphic to the universal enveloping
algebra of $\f$. Let us compare the algebra $\A/(\fr)$ (it is isomorphic
to the universal enveloping algebra of $\g$) and the algebra
$\B=\sum_{k=0}^\infty I^k/I^{k+1}$.

Since $R\in\gamma_2F$, it is readily seen that $\B_2=\A_2/(\fr)_2$ and
$\B_3=\A_3/(\fr)_3$. By the Poincar\'e-Birkhoff-Witt theorem (which is
valid over $\zz$, see \cite{L}), $\g$ is imbedded in $\A/(\fr)$.
Therefore, $\g_2$ is imbedded in $\B_2$ and $\g_3$ is imbedded in $\B_3$.
The proposition for $k=3$ and $4$ follows.
\end{proof}

Consider the algebra $\B=\sum_{k=0}^\infty I^k/I^{k+1}$. Note that
$\B_0=\zz$ and $\B_1=H_1$. Using Proposition \ref{digr}, it is easily
shown that $\B_2=H_1\otimes H_1/\bar\Delta(H_2)$ and $\B_3=H_1\otimes
H_1\otimes H_1/(\bar\Delta(H_2)\otimes H_1+H_1\otimes\bar\Delta(H_2))$.
Denote by $P_2$ (resp.\ $P_3$) the image of $[H_1,H_1]\subset H_1\otimes
H_1$ (resp.\ $[H_1,[H_1,H_1]]\subset H_1\otimes H_1\otimes H_1$) under
the canonical projection $H_1\otimes H_1\to H_1\otimes
H_1/\bar\Delta(H_2)$ (resp.\ $H_1\otimes H_1\otimes H_1\to H_1\otimes
H_1\otimes H_1/(\bar\Delta(H_2)\otimes H_1+H_1\otimes\bar\Delta(H_2))$).
By Proposition \ref{gamma}, the natural homomorphisms
$\gamma_2G/\gamma_3G\to \B_2$ and $\gamma_3G/\gamma_4G\to \B_2$ are
injective. Therefore, $\gamma_2G/\gamma_3G$ is isomorphic to $P_2$, and
$\gamma_3G/\gamma_4G$ is isomorphic to $P_3$.

Let $w$ be an element of $H_2$. We have
$\bar\Delta(w)=\sum_{i<j}\alpha_{ij}[h_i,h_j]$. Let
$$
\tau(w)=\prod_{i=1}^{n-1}\prod_{j=i+1}^n
(g_ig_jg_i^{-1}g_j^{-1})^{\alpha_{ij}}\in G.
$$
Clearly, $\tau(w)\in\gamma_3G$. Denote by $\bar\tau(w)$ the
corresponding element of $\gamma_3G/\gamma_4G$. Thus, we get a
homomorphism of Abelian groups $\bar\tau:H_2\to P_3$.

Note that $$R^2\simeq\Hom_\zz(H_1\otimes H_1/\bar\Delta(H_2),\zz)$$ and
$$Q^3=\Hom_\zz(H_1\otimes H_1\otimes H_1/(\bar\Delta(H_2)\otimes
H_1+H_1\otimes\bar\Delta(H_2)),\zz).$$ We have
$P_2\simeq[H_1,H_1]/{\Bar\Delta}(H_2)$. Let $\jmath:\Lambda^3H_1\to
P_2\otimes H_1$ be given by $\jmath(x\wedge y\wedge z)=[x,y]\otimes
z+[y,z]\otimes x+[z,x]\otimes y$. We have $P_3=(P_2\otimes
H_1)/\jmath(\Lambda^3H_1)$. Clearly, $\bar R^2=\Hom_\zz(P_2,\zz)$ and
$\bar Q^3=\Hom_\zz(P_3,\zz)$.

Since $X=B_.G$ is a simplicial set with single vertex, we see that the
map $\ka:H^1\to Z^1$ right inverse to the canonical projection is
unique. It is readily seen that there is a unique map $\imath:C_1\to
H_1\otimes H_1/\bar\Delta(H_2)$ such that $\imath(g_i)=0$ for
$i=1,\dots,n$ and the diagram
 $$
 \begin{CD}
 C_2 @>\Delta>> C_1\otimes C_1\\\
 @V{\dd}VV    @VVV\\
 C_1 @>\imath>> H_1\otimes H_1/\bar\Delta(H_2)
 \end{CD}
 $$
 is commutative. We set $\nu:R^2\to C^1$ to be the conjugate of
$\imath$.

Note that $\nu(\xi_i\otimes \xi_i)=\zeta(\ka(\xi_i))$ for $i=1,\dots,n$
and $\nu(\xi_i\otimes \xi_j+\xi_j\otimes
\xi_i)=\zeta(\ka(\xi_i+\xi_j))-\zeta(\ka(\xi_i))-\zeta(\ka(\xi_j))$ for
$i\ne j$ (where $\zeta(\omega)(a)=(\omega(a)-\omega(a)^2)/2$ for any
1-simplex $a$). Indeed, for any $r\in R^2$, the condition
$d\nu(r)=\mu\circ(\ka\otimes\ka)(r)$ means that $\nu(r)$ is determined
uniquely by its values on $1$-cycles $(g_1),\dots,(g_n)$. By the
definition of $\nu$, these values are zero for any $r\in R^2$. But
$\zeta(\ka(\xi_i))(g_j)=(\ka(\xi_i)(g_j)-\ka(\xi_i)(g_j)^2)/2=0$ for any
$i,j$, since 
$$\ka(\xi_i)(g_j)=
\begin{cases}
1,&\text{for }i=j,\\
0,&\text{for }i\ne j.
\end{cases}
$$
Similarly, $\zeta(\ka(\xi_i)+\ka(\xi_j))(g_k)=0$ for any $i,j,k$, $i\ne j$.

Thus, we can construct the map $\Bar{\Bar\la}:\bar Q^3\to H^2$ as above.

\begin{pr}
The map $\Bar{\Bar\la}:\bar Q^3\to H^2$ is conjugate
to the map $\bar\tau:H_2\to P_3$.
\end{pr}

\begin{proof}
Let $w\in H_2$; $\bar\Delta(w)=\sum_{i<j}\alpha_{ij}[h_i,h_j]$.
We want to describe $\bar\tau(w)\in P_3\subset \B_3$ in terms of
comultiplication in the chain complex of $X$.

Note that $\B_3$ is isomorphic to the kernel of the projection
$D^{(4)}(G)\to D^{(3)}(G)$. 
We identify $D^{(4)}(G)$ with $A^{(4)}(X)$ (by Proposition \ref{digr}).
Denote the image of $(g_i)\in X_1$ in
$A^{(4)}(G)$ by $a_i$. Then the image of $\tau(w)$ in $A^{(4)}(G)$ is
equal to
\begin{multline*}
\prod_{i<j}(1+a_i)(1+a_j)(1-a_i+a_i^2-a_i^3)(1-a_j+a_j^2-a_j^3)\\=
1+\sum_{i<j}\alpha_{ij}(a_ia_j-a_ja_i)(1-a_i-a_j).
\end{multline*}
Let us consider $\bar\tau$ as the map $H_2\to\B_3$. We see that it is
the sum of two maps---$\bar\tau_1$ and $\bar\tau_2$---where
$$
\bar\tau_1(w)=-\sum_{i<j}\alpha_{ij}(h_i\otimes h_j-h_j\otimes
h_i)\otimes(h_i+h_j)+(\bar\Delta(H_2)\otimes
H_1+H_1\otimes\bar\Delta(H_2))
$$
and $\bar\tau_2(w)$ is the preimage of
$$
\sum_{i<j}\alpha_{ij}(a_ia_j-a_ja_i)\in A^{(4)}(G)
$$
under the injection $\B_3\to A^{(4)}(G)$.

The element $\bar\tau_2(w)$ can be described as follows.  Let $c\in Z_2$
be a representative of the class $w$. The image of $\Delta c\in
C_1\otimes C_1$ in $A^{(4)}(X)$ equals zero. Therefore, in the
definition of $\bar\tau_2(w)$ we can use the image of  the element
$t=-\Delta
c+\sum_{i<j}\alpha_{ij}((g_i)\otimes(g_j)-(g_j)\otimes(g_i))\in
C_1\otimes C_1$ in $A^{(4)}(X)$ instead of
$\sum_{i<j}\alpha_{ij}(a_ia_j-a_ja_i)$ (which is the image of
$\sum_{i<j}\alpha_{ij}((g_i)\otimes(g_j)-(g_j)\otimes(g_i))\in
C_1\otimes C_1$). Clearly,  $t$ is the projection of $-\Delta c$ to
$C_1\otimes B_1+B_1\otimes C_1$ along the linear span of
$(g_i)\otimes(g_j)$. Hence, there is an element $\tilde t\in C_2\otimes
C_1\oplus C_1\otimes C_2$ such that
$t=(\dd\otimes\id\oplus\id\otimes\dd)\tilde t$. From the definition of
$A^{(4)}(X)$ it follows that $\bar\tau_2(w)$ is the image of
$-(\Delta\otimes\id\oplus\id\otimes\Delta)\tilde t$ under the natural
projection $C_1\otimes C_1\otimes C_1\to \B_3$. Consequently,
$\bar\tau_2(w)=\vphi(\imath\otimes pr+pr\otimes\imath)\Delta
c$, where $pr$ is the canonical projection $C_1\to H_1$ and $\vphi$ is
the natural map $ (H_1\otimes H_1/\bar\Delta(H_2))\otimes H_1\oplus
 H_1\otimes(H_1\otimes H_1/\bar\Delta(H_2))\to\B_3$.

Thus, we see that
$\bar\la$ is the conjugate of $\bar\tau_2$. Note that $\bar\tau_1(w)$ is
orthogonal to $p\bar Q^3$. On the other hand, the injections $p:\bar
Q^3\to Q^3$ and $P_3\to\B_3$ preserve the natural pairing. Since the
image of $\bar\tau=\bar\tau_1+\bar\tau_2$ belongs to $P_3$, the theorem
follows.
\end{proof}

\begin{thm}
Suppose that the the groups $G_a=\pi_1(U_a,u_a)$ and
$G_b=\pi_1(U_b,u_b)$ satisfy conditions (1)--(3). We identify
$H_1(U_a)$ with $G_a/\gamma_2G_a$ and $H_1(U_b)$ with $G_b/\gamma_2G_b$.
Let $f:H_1(U_a)\to H_1(U_b)$ be any isomorphism and let
$\vphi:H^1(U_b)\to H^1(U_a)$ be its conjugate. Then we claim the
following.
\begin{enumerate}
\item
The isomorphism $f$ can be extended to an isomorphism
$$G_a/\gamma_3G_a\to G_b/\gamma_3G_b$$
if and only if $\Lambda^2\vphi(\bar R^2_b)=\bar R^2_a$.
\item
Suppose that the previous condition hold and $P_3^{(a)}\cong P_3^{(b)}$
is a free Abelian group. Then the isomorphism $f$ can be extended to an
isomorphism
$$G_a/\gamma_4G_a\to G_b/\gamma_4G_b$$
if and only if
$[\la_a]$ corresponds to $[\la_b]$, i.~e., there is a commutative
diagram\footnote{This diagram is typeset with the help of \Xy-pic
package.}
$$
  \xymatrix{
 \bar Q^3_a
 \ar[d]^{\Bar{\Bar\lambda}_a}
 \ar[dr]
& & 
\bar Q^3_b
 \ar[ll]_-{\vphi\otimes\Lambda^2\vphi}
 \ar[d]^{\Bar{\Bar\lambda}_b}
 \ar[dl]
\\
 H^2_a
 & H
   \ar[l]
   \ar[r]
 &
 H^2_b
  }
$$
\end{enumerate}
\end{thm}

\begin{proof}
Clear.
\end{proof}

\begin{rem}
It is readily seen that $[\lambda]$ is just the invariant used in
\cite{Ry} to distinguish fundamental groups of combinatorially
equivalent complex hyperplane  arrangements.
\end{rem}


\begin{thebibliography}{88}

\bibitem {A} J.~F.~Adams, {\it On the cobar construction},
Colloque de topologie algebraique (Louvain,1956), George Thone, 1957,
pp.~81--87.

\bibitem{C1} K.-T.~Chen, {\it Iterated integrals, fundamental groups,
and covering spaces\/}, Trans.\ Amer.\ Math.\ Soc., {\bf 206} (1975),
pp. 83--98.

\bibitem{C2} K.-T.~Chen, {\it Extensions of $C^\infty$ function algebra
by integrals and Malcev completion of $\pi_1$\/}, Advances in Math.,
{\bf 23} (1977), pp. 181--210.

\bibitem{GM} P.~A.~Griffiths and J.~W.~Morgan, {\it Rational Homotopy
Theory and Differential Forms}, Progress in Math, {\bf 16} (1981).

\bibitem{Ko} T.~Kohno, {\it On the holonomy Lie algebra and the
nilpotent completion of the fundamental group of the complement of
hypersurfaces}, Nagoya J. Math. {\bf 92} (1983), 21--37.

\bibitem {L} M.~Lazard, {\it Sur les alg\`ebres enveloppantes
universelles de certaines alg\`ebres de Lie},
C.\ R., {\bf 235} (1952), pp. 788--791.

\bibitem {LW} A. Lundell and S. Weingram, {\it The topology of CW complexes},
NY, 1969.

\bibitem {M} W. Magnus, {\it \"Uber Beziehungen zwischen h\"oheren
Kommutattoren}, J.~Reine Angew.\ Math., {\bf 177} (1937), pp. 105--115.

\bibitem{MKS} W.~Magnus, A.~Karrass and D.~Solitar,
{\it Combinatorial group theory}, Interscience, New York etc., 1966.

\bibitem {OT}  P.~Orlik and H.~Terao, {\it Arrangements of hyperplanes},
Springer-Verlag, Berlin etc., 1992.

\bibitem{Q} D.~Quillen, {\it Rational homotopy theory},
Ann.\ Math. (2), {\bf 90} (1969), pp. 205--295.

\bibitem{R} I.~A.~Rips, {\it On the fourth integer dimension subgroup},
Israel J. Math. {\bf 12} (1972), pp. 342--346.

\bibitem{Ry} G.~Rybnikov, {\it On the fundamental group of the
complement of a complex hyperplane arrangement},
DIMACS Tech. Report 94-13 (1994), pp. 33--50
(math.AG/9805056).

\bibitem {Sa} R.~Sandling, {\it The dimension subgroup problem},
J. Algebra {\bf 21} (1972), pp. 216--231.

\bibitem {S} V.~A.~Smirnov, {\it On the cochain complex of topological
spaces}, Matem.\ Sbornik (Russian) {\bf 115} (1981), pp. 146--158.

\bibitem {Su1} D.~Sullivan, {\it Topology of manifolds and differential
forms}, Proceedings of conference on manifolds, Tokio, 1973.

\bibitem {Su2} D.~Sullivan, {\it Infinitesimal computations in
topology}, Publications de IHES, {\bf 47} (1977), pp. 269--331.

\end{thebibliography}
\end{document}